\documentclass[a4paper,11pt]{amsart}
\usepackage{amssymb,amscd}
\newtheorem{lemma}{Lemma}
\newtheorem{prop}{Proposition}
\newtheorem{thm}{Theorem}

\theoremstyle{definition}

\newtheorem{defn}{Definition}
\theoremstyle{remark}
\newtheorem{rem}{Remark}

\newcommand{\R}{{\mathbb R}}
\newcommand{\C}{{\mathbb C}}
\newcommand{\Z}{{\mathbb Z}}

\newcommand{\g}{\mathfrak g}
\newcommand{\n}{\mathfrak n}
\renewcommand{\t}{\mathfrak t}

\newcommand{\restr}[1]{|_{#1}^{\vphantom x}}
\begin{document}
\title[The Guillemin formula]
{The Guillemin formula and K{\"a}hler metrics\\ on toric symplectic manifolds}
\author[D. Calderbank, L. David and P. Gauduchon]
{David M. J. Calderbank, Liana David and Paul Gauduchon}
\thanks
{The second author is supported by by the Leverhulme Trust
and the William Gordon Seggie Brown Trust. All three authors are members
of EDGE, Research Training Network HPRN-CT-2000-00101, supported by the
European Human Potential Programme.}
\address{David M. J. Calderbank \\ Department of Mathematics and
  Statistics\\ University of Edinburgh\\ King's Buildings\\ Mayfield
  Road\\ Edinburgh EH9 3JZ\\ Scotland}
\email{davidmjc@maths.ed.ac.uk}
\address{Liana David and Paul Gauduchon \\ Centre de Math{\'e}matiques\\
{\'E}cole Polytechnique \\ UMR 7640 du CNRS
\\ 91128 Palaiseau \\ France}
\email{pg@math.polytechnique.fr}
\curraddr{Liana David\\ Department of Mathematics\\ University of Hull\\
Cottingham Road\\ Hull HU6 7RX\\ UK}
\email{L.R.David@hull.ac.uk}
\date{\today}
\begin{abstract}
  We discuss the construction of toric K{\"a}hler metrics on symplectic
  $2n$-manifolds with a hamiltonian $n$-torus action and present a
  simple derivation of the Guillemin formula for a distinguished
  K{\"a}hler metric on any such manifold. The results also apply to orbifolds.
\end{abstract}
\maketitle

\section*{Introduction}

There is a one-to-one correspondence, established by T.~Delzant in \cite{D},
between the class of compact connected symplectic manifolds of dimension
$2n$ with an effective hamiltonian action of the $n$-dimensional torus $T
^n$, known as compact {\it toric} symplectic manifolds, and a class a convex
polytopes in ${\mathbb R} ^n$, known as {\it Delzant polytopes}. More
precisely, given a compact toric symplectic manifold, the image of the
momentum map is a Delzant polytope; conversely, to each Delzant polytope
$\Delta$ is {\it canonically} associated a compact toric symplectic manifold
$M _{\Delta}$, in such a way that $\Delta$ may be identified with the image
of the corresponding momentum map.  By extending the class of Delzant
polytopes, a similar correspondence can be established for
orbifolds~\cite{LT}.

It turns out that the symplectic manifold $M _{\Delta}$ canonically
associated to $\Delta$ comes naturally equipped with a K{\"a}hler structure
$(M _{\Delta}, g, J)$ and a holomorphic action of the complexified torus
$T^n _{\C}$, which makes the complex manifold $(M _{\Delta}, J)$ into a
toric variety; moreover, on the dense open orbit of $T ^n _{\C}$, the
K{\"a}hler structure admits a $T ^n$-invariant globally defined K{\"a}hler
potential, which has been computed by V.~Guillemin in \cite{VG1} (see also
\cite{VG2} and \cite{abreu}).

The aim of this note is to provide a simple derivation of the Guillemin
formula. The idea is to systematically exploit the symplectic geometry of
the Delzant construction, in which the symplectic manifold $M_\Delta$ is
obtained as a symplectic quotient of $\R^{2d}$ (equipped with its standard
hamiltonian action of $T^d$) by a $(d-n)$-subtorus of $T^d$ (here $d$ is the
number of codimension one faces of the Delzant polytope).  The canonical
K{\"a}hler metric is obtained from the identification of $\R^{2d}$ with
$\C^d$ induced by the action of $T^d$: the flat metric on $\C^d$ then
descends to a K{\"a}hler metric on $M_\Delta$, the Guillemin metric.

The key observation is that a K{\"a}hler metric on a toric symplectic
manifold determines and---modulo the choice of angular coordinates---is
determined by a {\it reduced} metric on the interior of the image of its
momentum map, and that these reduced metrics behave in a straightforward way
with respect to symplectic quotients. In particular, in the Delzant
construction, the dual $\R^{n*}$ of the Lie algebra of the quotient torus
$T^n$ acting on $M_\Delta$ is naturally an affine subspace of the dual
$\R^{d*}$ of the Lie algebra of $T^d$, and the reduced metric on the
interior $\Delta_0$ of $\Delta$ is simply the pullback of the reduced metric
induced (on a quadrant of $\R^{d*}$) by the flat metric on $\C^d$. To carry
out this pullback, one simply has to write the flat metric on $\C^d$ in {\it
momentum coordinates}, i.e.,
\begin{equation}
g_0 = \sum_{j=1}^d \frac{{d\mu_j}^2}{2\mu_j} + \sum_{j=1}^d 2\mu_j dt_j^2.
\end{equation}
(Here the first sum gives the reduced metric, the second the metric on the
$d$-torus fibres.)

This procedure leads naturally to the Guillemin metric on $M_\Delta$ in
momentum coordinates, which are the derivatives of a K{\"a}hler potential
with respect to the more traditional {\it holomorphic coordinates}
associated to the complexified torus action. The Legendre transform $G$ of
the K{\"a}hler potential $F$ gives a {\it dual} or {\it symplectic}
potential for the toric K{\"a}hler metric in momentum
coordinates~\cite{VG1,VG2}.  The usefulness of dual potentials and momentum
coordinates in toric symplectic geometry has been emphasised by
M.~Abreu~\cite{abreu,abreu2}. In particular, dual potentials, like reduced
metrics, behave straightforwardly under symplectic quotients, and so the
dual potential of $M_\Delta$ is easy to compute.  Conversely, inverting the
Legendre transform, we can reobtain the holomorphic coordinates and the
K{\"a}hler potential $F$, yielding the Guillemin formula.

Although this indirect approach is very simple, it seems not to have been
described in detail before. However, it is closely related to the
construction of toric hyperk{\"a}hler metrics in \cite{BD} by R.~Bielawski and
A.~Dancer, using the method of~\cite[Section 2(C)]{HKLR}. As they point out,
the hyperk{\"a}hler approach also gives the Guillemin formula by setting to
zero half of the coordinates. Nonetheless, we hope it will be of use to
present directly the K{\"a}hler story, with complete arguments, here.

We begin by reviewing the Delzant--Lerman--Tolman theory of toric symplectic
manifolds and orbifolds. Then we introduce toric K{\"a}hler metrics and
their description in momentum coordinates. We prove that under symplectic
quotient, the reduced metrics are related by pullback, then apply this to
the flat metric on $\C^d$.

Finally, in an appendix, we explain how the Guillemin formula for the
K{\"a}hler potential is also a direct consequence of a general formula
appearing in \cite{BG} (cf.~also \cite[Section 3 (E)]{HKLR} ).

\section{Hamiltonian group actions} \label{s:ham}

Let $(M,\omega)$ be a symplectic manifold with an effective action of a Lie
group $G$. For simplicity, we assume that $G$ is compact.  Then we say the
action is hamiltonian if there is a momentum map $\mu\colon M\to\g^*$ for
the action, i.e., $\mu$ is $G$-equivariant and $\iota_{X_\xi}\omega=-\langle
d\mu,\xi\rangle$ for any $\xi\in\g$, where $X_\xi$ is the corresponding
vector field on $M$ (generating the action of a $1$-parameter subgroup of
$G$).

If $c$ is a regular value of $\mu$ invariant under the coadjoint action of
$G$, then $\mu^{-1}(c)$ is a $G$-invariant submanifold of $M$ and the action
of $G$ on $\mu^{-1}(c)$ is locally free, with finite isotropy groups. If
this action is free, then $M{/\!/\!}_c G:=\mu^{-1}(c)/G$ is a symplectic
manifold, the {\it symplectic quotient} of $M$ by $G$. In general, however,
$M{/\!/\!}_cG$ is a symplectic orbifold.

Symplectic reduction can also be used to construct symplectic manifolds and
orbifolds with hamiltonian group actions, by taking a symplectic quotient by
a normal subgroup of the symmetry group.

\begin{lemma} \label{partquot}
Let $G$ be a (compact) Lie group with closed Lie subgroup $N$, and suppose
that $(\tilde M,\tilde \omega)$ has a hamiltonian $G$-action with momentum
map $\mu\colon\tilde M\to\g^*$, so that $\mu_N=i^*\circ\mu\colon \tilde M\to
\n^*$ is the momentum map for the induced hamiltonian $N$-action (here $i^*$
is adjoint to the inclusion of $\n$ in $\g$). Suppose further that $N$ is
normal in $G$, that $\tilde c\in \g^*$ is invariant under the coadjoint
action of $G$, and that $c=i^*(\tilde c)$ is a regular value of $\mu_N$ such
that the action of $N$ on $\mu_N ^{-1}(c)$ is locally free.

Then $M=\tilde M{/\!/\!}_c N$ has a hamiltonian action of $G/N$ whose
momentum map, viewed as an $N$-invariant function on $\mu _N ^{-1} (c)$, is
the restriction of $\mu - \tilde{c}$ to $\mu _N ^{-1} (c)$.
\end{lemma} 
\begin{proof} The action of $[g]\in G/N$ on $[z]\in M=\mu _N ^{-1} (c)/N$
is given by $[g] \cdot [z] = [g \cdot z]$. To see that this is well defined,
note that the elements $z$ of $\mu _N ^{-1} (c)$ are characterized by
$\langle \mu (z) - \tilde{c}, \zeta \rangle = 0$ for all
$\zeta\in \n$; now since $\tilde{c}$ is $G$-invariant, for any $z$ in $\mu
_N ^{-1} (c)$ and any $g$ in $G$, we then have $\langle \mu (g \cdot z) -
\tilde{c}, \zeta \rangle = \langle \mu (z) - \tilde{c}, {\rm Ad} _{g ^{-1}}
\zeta \rangle$. Thus, since $N$ is normal in $G$, $g\cdot z$ is in $\mu _N
^{-1} (c)$ and $[g\cdot z]$ depends only on $[g]$ and $[z]$.

The fact that $\langle \mu (z) - \tilde{c}, \zeta \rangle = 0$ for any $z\in
\mu _N ^{-1} (c)$ means that the restriction of $\mu - \tilde{c}$ to $\mu _N
^{-1} (c)$, say $\nu$, has its values in $(\g/\n)^*$, identified with the
annihilator $\n ^0$ of $\n$ in $\g^*$.  Moreover, $\nu$ is $N$-invariant:
since $N$ is normal, ${\rm Ad}^*_g$ acts by the identity on $\n^0$ for any
$g$ in $N$.

Now $\nu$ gives a momentum map for the action of $G/N$ on $M$, since for any
$\xi\in\g$, $\langle d \nu, [\xi] \rangle$ pulls back to $\langle d\mu, \xi
\rangle = -\iota_{X_\xi}\tilde \omega$ on $\mu _N ^{-1} (c)$.
\end{proof}

We remark that if $u^*\colon (\g/\n)^*\to \g^*$ is the natural inclusion
(adjoint to $u\colon \g\to\g/\n$) then the affine embedding $\ell=u^*+\tilde
c$ gives a commutative diagram
\begin{equation}
\begin{CD} 
  \mu _N ^{-1} (c) @> \mu >> (i^*)^{-1}(c) \,\hbox to 0pt{$\subset \g^*$} \\
  @V{q}VV         @AA{\ell }A  \\
  M                @>\nu >>  (\g/\n)^{*\!\!}
\end{CD}
\end{equation}
where $q$ is the natural projection. In particular $\ell$ identifies the
image of $\nu$ with the intersection of the affine subspace $(i^*)^{-1}(c)$
and the image of $\mu\colon M\to\g^*$.

\section{Toric symplectic manifolds} \label{s:toric}

A {\it hamiltonian $m$-torus action} on a symplectic manifold $(M,\omega)$
of dimension $2n$ is an effective torus action generated by an
$m$-dimensional family of hamiltonian vector fields $X \in
C^\infty(M,TM)\otimes\R^{m*}$ with a momentum map $\mu\colon M\to \R^{m*}$
whose components Poisson commute, i.e., $\iota_X\omega = -d\mu$ and
$\omega(X,X)=0$. (This ensures that the components of $X$ commute, and that
$\mu$ is $T^m$-equivariant.)

Observe that $\R^m$ is here the Lie algebra of the torus, and so contains a
lattice $\Z^m$ such that the torus is $T^m=\R^m/2\pi\Z^m$.

Since the orbits are isotropic with respect to $\omega$, it follows that
$m\leq n$---if equality holds, we say that $(M,\omega,\mu)$ is a {\it toric}
symplectic manifold.

In this section we summarize the classification of {\it compact} toric
symplectic manifolds, established by Delzant~\cite{D}, and the extension of
this theory to orbifolds~\cite{LT}. It turns out that there is a one to one
correspondence between compact connected toric symplectic manifolds (or
orbifolds) and a class of compact convex polytopes, called {\it Delzant
polytopes}. More precisely, given a compact connected toric symplectic
manifold $(M^{2n},\omega,\mu)$, the image of the momentum map $\mu$ is a
compact convex polytope $\Delta$ in $\R^{n*}$ by the
Atiyah--Guillemin--Sternberg convexity theorem. The same holds for toric
orbifolds; the polytopes obtained belong to the following class.

\begin{defn} Let $\R^n$ be an $n$-dimensional vector space with a lattice
$\Z^n$, and consider a compact convex polytope $\Delta$ in $\R^{n*}$
defined by the equations
\begin{equation} \label{polytope}
\langle x, u_j\rangle  \geq \lambda_j,
\end{equation}
($u_j\in\R^n$, $\lambda_j\in\R$) for $j=1,\ldots d$, where $d > n$ is the
number of $(n - 1)$-faces (simply called {\it faces} in the sequel) of
$\Delta$.

Then $\Delta$ is said to be a {\it rational Delzant polytope} if the $u_j$
belong to the lattice $\Z^n$ and the elements corresponding to faces
meeting any given vertex of the polytope form a basis of $\R^n$. (It
follows that the polytope is $n$-valent, i.e., $n$ edges meet each
vertex---the dual basis of $\R^{n*}$ gives the directions of these edges.)

The polytope is said to be {\it integral} or simply a {\it Delzant polytope}
if the $u_j$ corresponding to faces meeting any vertex form a ($\Z$-)basis
of the lattice $\Z^n$.
\end{defn}

\begin{rem}
Note that in this definition we consider the vectors $u_j$ and the lattice
$\Z^n$ as part of the data. If $\Delta$ is integral, then in fact the $u_j$
are uniquely determined by $\Delta$ and the lattice, since they must be
primitive; conversely the $u_j$ generate the lattice in this case.

More generally, the $u_j$ are determined by $\Delta$, $\Z^n$ and a positive
integer labelling of the faces of $\Delta$, since we may write $u_j=m_j y_j$
with $y_j$ primitive and $m_j\in\Z^+$---because of this description, rational
Delzant polytopes are called {\it labelled polytopes} in~\cite{LT}. Now,
conversely, the $u_j$ in general only generate a (finite index) sublattice
$\hat\Z^n$ of $\Z^n$. However, there is not much generality lost by replacing
$\Z^n$ with this sublattice---see Remark~\ref{orbifold} below.
\end{rem}

The image of the momentum map of a symplectic toric orbifold (manifold) is a
rational (integral) Delzant polytope. Conversely, to each rational Delzant
polytope is canonically associated a compact connected toric symplectic
orbifold $(M _{\Delta},\omega_\Delta,\mu_\Delta)$ such that $\Delta$ is the
image of the momentum map $\mu_\Delta$.  A toric symplectic orbifold is
uniquely determined by its polytope up to equivariant symplectomorphism, so
the two operations are mutually inverse.

The construction of $M _{\Delta}$ will be important in the sequel, so we
summarize it here, following~\cite{VG2,LT}. The idea is to obtain $M_\Delta$
as a symplectic quotient of a $2d$-dimensional symplectic vector space by
the natural action of a $(d-n)$-dimensional subgroup $N$ of
$T^d=\R^d/2\pi\Z^d$.

Here, $\Z^d$ may be thought of (canonically) as the free abelian group
generated by the $d$ faces $\sigma_j$ of $\Delta$, and $\R^d$ is the
corresponding real vector space $\Z^d\otimes_\Z\R$. We take the symplectic
vector space to be $\R^{2d}=\Z^d\otimes_\Z\R^2$, where $\R^2$ is equipped
with its standard area form and circle action. In terms of polar coordinates
on $\R^2\cong\C$, we then have a natural action of the torus $T^d$ on
$\R^{2d}$:
\begin{equation}\label{tact}
[a_1, \ldots, a _d] \cdot \sum_{j=1}^d \sigma_j\otimes(r_j,\theta_j) =
\sum_{j=1}^d \sigma_j\otimes(r_j,\theta_j+a_j)
\end{equation}
where $[a_1, \ldots, a _d]$ denotes the class of $(a_1, \ldots, a _d) \in
\R^d$ modulo $2\pi\Z^d$. This action is hamiltonian with respect to the
standard symplectic form, which in polar coordinates is $\sum_{j=1}^d r_j
dr_j\wedge d\theta_j$. It follows that the momentum map $\mu$ of the $T^d$
action is given (up to translation) by
\begin{equation} \label{J0}
\mu(v) =\frac12 (r_1^2,\ldots r_d^2)=\frac12 \sum_{j=1}^d r_j^2 \sigma_j^*
\end{equation}
where $\sigma_j^*$ is the basis of $\R^{d*}$ dual to $\sigma_j$, and
$v=\sum_{j=1}^d \sigma_j\otimes(r_j,\theta_j)$.

Let $u\colon \Z ^d \to \Z ^n$ be the group homomorphism determined by
$\sigma _j \mapsto u _j$.  We denote also by $u$ the corresponding
homomorphisms from $\R^d$ to $\R^n$ and from $T^d$ to $T^n$, which are
surjective, and we denote the kernels by $\n$ and $N$ respectively.

\begin{rem} $u\colon \Z ^d \to \Z ^n$ is surjective if and only if
the $u_j$ generate $\Z^n$, in which case $N$ is the quotient of its
(abelian) Lie algebra $\n$ by the kernel of $u\colon 2\pi\Z ^d \to
2\pi\Z ^n$, and so is a subtorus of $T^d$.

More generally if the $u_j$ generate a sublattice $\hat\Z^n$, then the
quotient of $N$ by the connected component of the identity is the
finite group $\Z^n/\hat\Z^n$.
\end{rem}

We now consider the restriction of the action of $T^d$ on $\R^{2d}$ to the
subgroup $N$. The action of $N$ is still hamiltonian with momentum map
$\mu_N = i^* \circ \mu$, where $i^*\colon\R^{d*}\to \n^*$ is adjoint to the
inclusion $i\colon\n\to\R^d$.  We let $c= i^*(-\lambda)$, where $\lambda =
(\lambda _1,\ldots \lambda _d)$ is naturally the element $\sum_{j=1}^d
\lambda_j \sigma_j^*$ of $\R^{d*}$, and consider the momentum level set
$\mu_N^{-1}(c)$.

The annihilator in $\R ^{d*}$ of $\n$ is the image of the adjoint $u^*$ of
$u$.  It follows that $\mu_N ^{-1} (c)$ is the set of elements $v$ of $\R
^{2d}$ for which there exists $x$ in $\R^n$ with $\mu (v) + \lambda = u ^*
(x)$.

In other words, since $u^* (x) = \sum_{j=1}^d \langle x, u _j\rangle \,
\sigma _j^*$, $v=\sum_{j=1}^d \sigma_j\otimes(r_j,\theta_j)$
belongs to
$\mu_N ^{-1} (c)$ if and only if there exists $x$ in $\R^{n*}$ such that
\begin{equation} \label{mu}
\ell _j (x):= \langle x, u _j\rangle  - \lambda _j = \frac{r_j^2}{2},
\end{equation}
for $j = 1, \ldots, d$. 

We observe that $(\ell_1,\ldots \ell_d)$ are the components of the affine
embedding $\ell = u^*-\lambda$ of $\R^{n*}$ into $\R^{d*}$. Hence $v$
belongs to $\mu_N^{-1}(c)$ if and only if $\mu(v)$ belongs to the image of
$\ell$, and then $x$ appearing in~\eqref{mu} is uniquely determined by $v$
and belongs to $\Delta$.  This defines a map from $\mu_N ^{-1} (c)$ to
$\Delta$, which, following Lemma~\ref{partquot}, we denote by $\nu$. We thus
have:
\begin{equation} \label{phi}
\ell (\nu (v))=\mu(v)
\end{equation}
for all $v = \sum \sigma_j\otimes (r_j,\theta_j)$ in $\mu _N ^{-1} (c)$;
$\nu$ is clearly surjective and for each $x$ in $\Delta$ and $\nu ^{-1}
(x)$ coincides with the orbit of any of its elements under the action of
$T^d$. In particular, $\mu _N ^{-1} (c)$ is compact.

Notice that $\nu (v)$ belongs to the interior of $\Delta$ (the
complement of the union of its faces) if and only if all the $r _j$ are
non-zero; the isotropy group of $v$ in $T ^d$ is then zero.

It turns out that for {\it any} $v$ in $\mu_N ^{-1} (c)$, the isotropy group
$N _v = T _v \cap N$ of $v$ in $N$ is finite, and is zero if the polytope is
integral. Indeed, if $\nu (v)$ belongs to the intersection of $k$ faces
$\sigma_{j_1}, \ldots \sigma_{j_k}$, then $r _{j _1} = \ldots = r _{j _k} =
0$ and the isotropy group $T_v$ of $v$ in $T^d$ is the $k$-dimensional
subtorus of $T ^d$ whose Lie algebra $\t_v$ is the set of $(a_1,\ldots a_d)$
with $a_j=0$ for $j \neq j_1, \ldots, j _k$.  The intersection of $T_v$ with
$N$ is thus the set of $[a_1,\ldots a_d]\in T_v$ with $\sum _{j = 1} ^d a _j
u _j\in 2\pi\Z^d$.  However, the $u_j$ with $a_j\neq 0$ mod $2\pi\Z$, are
part of a basis for a (finite index) sublattice of $\Z^d$ (choose a vertex
in $\sigma_{j_1}\cap \ldots \cap \sigma_{j_k}$). This means that the $u_j$
with $a_j\neq 0$ mod $2\pi \Z$ form a basis for a finite index sublattice of
$\t_v\cap\Z^d$; $N _v$ is the quotient of $\t_v\cap\Z^d$ by this sublattice.

Therefore $N$ acts locally freely on $\mu_N ^{-1} (c)$, so $c$ is a regular
value of $\mu _N$, $\mu _N ^{-1} (c)$ is a closed submanifold of $\R ^{2d}$
and the quotient $M := \mu _N ^{-1} (c)/N$ is a compact symplectic orbifold,
of (real) dimension $2d - 2(d-n) = 2n$. This is the symplectic quotient, at
momentum level $c$, of $\R^{2d}$ by $N$.

By Lemma~\ref{partquot}, $\nu$ induces a momentum map with image $\Delta$
for the natural action of the quotient torus $T^n=T^d / N$ on $M = \mu _N
^{-1}(c)/N$.

\begin{rem} \label{orbifold}
If the $u_j$ do not generate $\Z^n$, then we can instead work with the
sublattice $\hat\Z^n$ generated by the $u_j$: $\hat N$ is then the connected
component of the identity in $N$, so $\hat M$ is a finite cover of $M$, with
a hamiltonian action of the cover $\hat T^n=\R^n/2\pi\hat\Z^n$ of
$T^n=\R^n/2\pi\Z^n$.  Thus $M$ is the quotient of $\hat M$ by the finite
group $\Z^n/\hat\Z^n$.  It is usual in the theory of orbifolds to pass to
such global covers if they exist, as this will (at least partially) resolve
some of the orbifold singularities. This explains why there is not much
generality lost in assuming that the $u_j$ generate $\Z^n$.
\end{rem}

\section{K{\"a}hler quotients}\label{s:quot}

In general, let $(M,g,J,\omega)$ be a K{\"a}hler manifold with an effective
action of a group $G$ by isometries, which is hamiltonian with momentum map
$\mu\colon M\to\g^*$. Let $M{/\!/\!}_c G$ be a symplectic quotient of $M$ by
$G$; then since $G$ acts by isometries on $\mu^{-1}(c)$ (and we are
supposing the action is locally free), $M{/\!/\!}_c G$ is a K{\"a}hler
orbifold, the {\it K{\"a}hler quotient} of $M$ by $G$.

We now consider the situation described by Lemma \ref{partquot}, where $G$
acts on $(\tilde M,\tilde\omega)$ and $N$ is a normal subgroup of $G$.  We
assume additionally that $G$ acts freely on $\tilde{M}$ and that the
dimension of $G$ is equal to half the dimension of $\tilde{M}$ so that
$\mu\colon\tilde{M}\to\g^*$ is a $G$-principal bundle over its image. Then
if $\tilde g$ is a $G$-invariant compatible metric on $\tilde M$, we obtain
a uniquely defined riemannian metric, say $\tilde{g} _{\rm red}$ on ${\rm
  im}\,\mu$ such that $\mu$ is a riemannian submersion from $\tilde{g}$ to
$\tilde{g} _{\rm red}$.

Similarly, the quotient group $G/N$ acts freely on the K{\"a}hler quotient
$M=\tilde M{/\!/\!}_c N$ and this makes $\nu\colon M\to (\g/n)^*$ into a
$G/N$-principal bundle over its image, and $\nu$ is a riemannian submersion
for a uniquely defined riemannian metric $g _{\rm red}$ on ${\rm im}\,\nu$,
when $M$ is equipped with the metric $g$ induced by $\tilde{g}$.

We denote by $\ell=u^*+\tilde c$ the affine map from ${\rm im}\,\nu$ to
${\rm im}\,\mu$ induced by the natural inclusion  $u^*\colon (\g/\n) ^*\to
\g^*$.

\begin{lemma} \label{redmetric}
The reduced metrics are related by $\ell ^* \tilde{g} _{\rm red} = g _{\rm
    red}$.
\end{lemma} 
\begin{proof} Let $z\in \mu _N ^{-1} (c)$
and let $\tilde X\in T_z\tilde M$ be orthogonal to the $G$-orbit through
$z$. Then the (isometric) projection $X$ of $\tilde X$ onto $T_x M$ is
orthogonal to the $G/N$-orbit through $x=q(z)$. We thus have
\begin{equation*}
|\nu _*(X)| _{g _{\rm red}} = |X| _g = |\tilde{X}|  _{\tilde{g}} =
|\mu_* (\tilde X))| _{\tilde{g} _{\rm red}}.
\end{equation*}
Since $\mu(z)=\ell(\nu(x))$ by Lemma \ref{partquot}, this is what we want.
\end{proof}

We wish to emphasise at this point that the K{\"a}hler quotient is a
construction in symplectic geometry---the riemannian metric on $M$ (and
hence the complex structure) is just a passenger; that is, {\it any}
$G$-invariant K{\"a}hler metric on $M$ with K{\"a}hler form $\omega$
descends to a compatible K{\"a}hler metric on $M{/\!/\!}_c G$.

In general it is difficult to parameterize $G$-invariant K{\"a}hler metrics
with a fixed K{\"a}hler form; instead one usually fixes the complex
structure, so that compatible K{\"a}hler metrics are parameterized (locally)
by {\it K{\"a}hler potentials} $K$ (i.e., $dd^cK=\omega$, where $d^c=J\circ
d$).  In the case of a compact toric symplectic manifolds however,
compatible toric K{\"a}hler metrics on the open set $M _0 = \mu ^{-1}
(\Delta _0)$ determine and are determined by a $T ^m$-invariant lagrangian
foliation (the integral leaves of the distribution generated by vector
fields $J K _1, \ldots, JK _m$) and by the induced riemannian metric on
$\Delta _0$, which is of hessian type, i.e. defined in terms of a {\it dual
potential}, which is related to the K{\"a}hler potential by Legendre
transform (more details in the next section).

\section{Dual potentials for toric K{\"a}hler metrics}\label{s:dual}

Suppose that $(M^{2n},\omega)$ is a toric symplectic manifold with momentum
map $x\colon M\to \R^{n*}$ and that $(g,J)$ is a compatible toric K{\"a}hler
structure. Let $X\in C^\infty(M,TM)\otimes\R^{n*}$ be the corresponding family
of hamiltonian Killing vector fields with $\iota_X\omega=-dx$.

We wish to study the geometry of $M$ locally, so we assume that the torus
action is free. Then $JX$ is a family of holomorphic vector fields and the
components of $X$ and $JX$ together form a frame of commuting vector fields.
The closed $1$-forms in the dual frame are naturally the components of
$\R^n$-valued $1$-forms $\alpha$ and $J\alpha=-\alpha\circ J$ with
$\alpha(X)$ equal to the identity of $\R^n$.  If we locally write
$\alpha=dt$ and $J\alpha=-dy$ with $y,t\colon M\to\R^n$, then $y+it\colon
M\to \C^n$ is a local holomorphic chart on $M$.

The components of $J\alpha$ and $\iota_X\omega$ span the same rank
$n$-subbundle of $T^*M$, the annihilator of tangent space to the orbits of
the torus action.  Hence there must be mutually inverse functions
$\boldsymbol F\colon M\to{\rm Hom}(\R^{n},\R^{n*})$ and $\boldsymbol G\colon
M\to{\rm Hom}(\R^{n*},\R^{n})$ such that
\begin{equation*}
\iota_X\omega=\langle \boldsymbol F,J\alpha\rangle\qquad\text{and}\qquad
J\alpha= \langle \boldsymbol G,\iota_X\omega\rangle.
\end{equation*}
Since $g(X,X)=\omega(X,JX)=\boldsymbol F$, we deduce that $\boldsymbol F$
and $\boldsymbol G$ are, at each point of $M$, symmetric and positive
definite elements of $\R^{n*}\otimes\R^{n*}$ and $\R^n\otimes\R^n$.

Concretely, with respect to a basis of the Lie algebra $\R^n$, we have
\begin{equation*}
dx_r= \sum_s F_{rs}dy_s\qquad\text{and}\qquad
dy_r= \sum_s G_{rs}dx_s.
\end{equation*}
Now $\sum_s G_{rs} dx_s$ is closed if and only if $G_{rs}$ is the
hessian with respect to $x$ of a function $G$. Indeed $\sum_s G_{rs}
dx_s=dy_r$ if and only if $G_{rs}=\partial y_r/\partial x_s$, in which
case $y_r=\partial G/\partial x_r$ for some function $G$, since $G_{rs}$
is symmetric. Similarly $\sum_s F_{rs} dy_s$ is closed if and only if
$F_{rs}$ is the hessian with respect to $y$ of a function $F$.

We therefore have functions $F,G$ on $M$ such that $x_r=\partial F/\partial
y_r$ and $y_r=\partial G/\partial x_r$; without loss of generality, we may
also assume that the constant $F+G-\langle x,y\rangle$ is zero. Thus the
coordinate systems $x$ and $y$, and the functions $F$ and $G$, are related
by Legendre transform:
\begin{equation*}
F+G = \langle x,y\rangle = \sum_r x_r \frac{\partial G}{\partial x_r}=
\sum_r y_r\frac{\partial F}{\partial y_r}.
\end{equation*}
$F$ is a K{\"a}hler potential, since $\omega=  \sum_{r,s} F_{rs} dy_s\wedge
dt_r=dd^cF$. On the other hand $G$ is the dual potential we need to
parameterize toric K{\"a}hler metrics with fixed symplectic form.

\begin{prop}\label{toric}~\textup{\cite{VG1,abreu}} Let $G$ be a function
on $\R^{n*}$ whose hessian $G_{rs}$ is positive definite with inverse
$F_{rs}$.  Then
\begin{equation*}
\sum_{r,s} (G_{rs} dx_r dx_s + F_{rs} dt_r dt_s)
\end{equation*}
is a toric K{\"a}hler metric with K{\"a}hler form
\begin{equation*}
\omega =  \sum_r dx_r\wedge dt_r
\end{equation*}
and momentum map $x$.
\end{prop}

Observe that $g_{\rm red}$ is the metric $\sum_{r,s} G_{rs} dx_r dx_s$ on the
image of the momentum map $x$. The following simple observation explains why
dual potentials are convenient when taking K{\"a}hler quotients of toric
K{\"a}hler manifolds.

\begin{prop} \label{toricquot} Let $(\tilde{M},\tilde g)$ be a toric K{\"a}hler
manifold of real dimension $2d$. Denote by $\tilde{M} _0$ the open set
on which the $T^d$ action is free, so that the momentum map $\mu\colon
\tilde{M} _0\to U\subset \R^{d*}$ is a principal $T^d$-bundle.

Let $N$ be a subgroup of $T ^d$ of codimension $n$, with momentum map,
$\mu_N=i ^*\circ \mu$, where $i^*\colon \R^{d*}\to\n^*$ is the natural
projection. Choose $\tilde{c}\in\R^{n*}$ such that $c=i^*(\tilde c)$ is a
regular value of $\mu _N$ and the action of $N$ on $\mu _N ^{-1} (c)$ is
locally free. Let $M = \mu _N ^{-1} (c)/N$ be the K{\"a}hler quotient with
metric $g$ and $q\colon \mu _N ^{-1} (c)\to M$ the natural projection.

The action of $T ^d$ on $\tilde{M}$ preserves $\mu _N ^{-1} (c)$,
hence induces a hamiltonian action of the quotient torus $T ^n := T ^d/N$ on
$M$ with momentum map $\nu$ defined by
\begin{equation}
\nu \circ q = \mu - \tilde{c}.
\end{equation}
Let $\Delta_0$ denote the image of $\nu$ restricted to $M_0$, the quotient
of $\mu _N ^{-1} (c) \cap \tilde{M} _0$ by $N$.

Then the reduced metric $g_{\rm red}$ on $\Delta_0$ is the pullback by
$\ell$ of the reduced metric $\tilde g_{\rm red}$ on $U$ and so if $\tilde
G\colon U\to\R$ is a dual potential for $\tilde g$ on $\tilde M_0$, then
$G=\tilde G\circ\ell \colon \Delta_0\to \R$ is a dual potential for $g$ on
$M_0$.
\end{prop}
\begin{proof} Most of the statements in this proposition are immediate from
Lemma~\ref{partquot}, since the coadjoint actions are trivial and $N$ is
necessarily a normal subgroup of $T^d$.  It remains to establish the
conclusion that $G=\ell^*\tilde G$ is a dual potential for the quotient.
This follows from Lemma~\ref{redmetric}, since $\ell$ is an affine map with
respect to the flat affine connections $D$ and $\tilde D$ on $\Delta_0$
and $U$. Indeed, for the dual potentials we have
\begin{equation*}
D d (\ell^*\tilde G) = \ell^*\tilde D d\tilde G =\ell^*\tilde g_{{\rm red}}
=g_{{\rm red}}.
\end{equation*}
\end{proof}

\section{The Guillemin formula}\label{s:guillemin}

In section~\ref{s:toric} we explained how any toric symplectic orbifold $M$
could be obtained as a symplectic quotient of a symplectic vector space
$\R^{2d}$. In canonical terms $\R^{2d}=\Z^d\otimes_\Z \R^2$, where $\Z^d$
denotes the free abelian group generated by the faces $\sigma_1,\ldots
\sigma_d$ of the (rational) Delzant polytope $\Delta$, and $\R^2$ is the
standard symplectic vector space with a (linear) circle action. To
emphasize the symplectic nature of the construction we deliberately
suppressed the obvious identification of $\R^2$ with $\C$, and hence of
$\R^{2d}$ with $\C^d$.  In these terms the torus action~\eqref{tact} on
$\C^d$ is given by
\begin{equation}
[a_1, \ldots, a _d] \cdot (z_1,\ldots z_d)
=(e^{ia_1}z_1,\ldots e^{ia_1}z_d)
\end{equation}
where $[a_1, \ldots a _d]\in\R^d/2\pi\Z^d$ and $z_j=r_j e^{i\theta_j}$.

Thus $\R^{2d}$ is in fact equipped with a canonical flat K{\"a}hler metric
compatible with the symplectic form $\frac{i}{2}\sum dz_j\wedge d\overline
z_j$, and this induces a canonical K{\"a}hler metric on $M$, called the {\it
  Guillemin metric}.

In this section, we obtain an explicit expression, the {\it Guillemin
  formula}, for the reduced metric on the interior $\Delta_0$ of the
(rational) Delzant polytope $\Delta$. We also obtain formulae for the
K{\"a}hler potential and dual potential. These results are due to
Guillemin~\cite{VG1,VG2} (see also Abreu~\cite{abreu,abreu2}), but we
obtain a simple new proof using Proposition~\ref{toricquot}.

Recall that the toric symplectic manifold $(M,\omega,\nu)$ fits into the
commutative diagram
\begin{equation}
\begin{CD} 
\mu _N ^{-1} (c) @> \mu >> \R ^{d*}\\
@V{q}VV         @AA{\ell }A  \\
M @>\nu >> \Delta \subset \R ^{n*}
\end{CD} 
\end{equation}
where $\ell\colon\R^{n*}\to\R^{d*}$ is the affine map
\begin{equation}
\ell (x) = (\ell _1 (x), \ldots, \ell _d (x))
\end{equation} 
with
\begin{equation}
\ell_j(x)=\langle u_j,x\rangle -\lambda_j.
\end{equation}

Now we restrict attention to the open subset of $\C^d$ where the $T^d$
action is free, and the corresponding open subset $M_0$ of $M$. These
are principal $T^d$ and $T^n$ bundles over $U\subset\R^{d*}$ and
$\Delta_0\subset\R^{n*}$ respectively, where $U$ is the positive
quadrant of $\R^{d*}$.

It remains only to compute the reduced metric $\tilde g_{\rm red}$ and dual
potential $\tilde G$ for the flat K{\"a}hler structure on $\C^d$ and pull
these back by $\ell$. The K{\"a}hler potential can then be found by Legendre
transform.

The explicit derivation of $\tilde g_{\rm red}$ and $\tilde{G}$ is a
straightforward coordinate transformation: $z_j=r_je^{i\theta_j}$ and the
components of the momentum map $\mu$ are given by $\mu_j=r_j^2/2$. Thus
the flat metric on $\C^d$ is
\begin{equation}
\sum_{j=1}^d dz_j\,d\overline z_j = \sum_{j=1}^d
\bigl({dr_j} ^2+r_j^2\, {d\theta_j}^2\bigr)
= \sum_{j=1}^d
\Bigl(\frac{{d\mu_j}^2}{2\mu_j}+2 \mu_j \,{d\theta_j}^2\Bigr).
\end{equation}
Hence
\begin{equation}\label{tgred}
\tilde g_{\rm red} = \frac12 \sum_{j=1}^d \frac{{d\mu_j}^2}{\mu_j}
\end{equation}
so that a dual potential $\tilde G$ is given by
\begin{equation}\label{tG}
\tilde G = \frac 12 \sum_{j=1}^d \mu_j\log \mu_j.
\end{equation}

\begin{thm} \label{thguillemin} Let $(M^{2n},\omega)$ be a toric symplectic
orbifold with momentum map $\nu\colon M\to \Delta\subset \R^{n*}$, where the
(rational) Delzant polytope $\Delta$ is given by
\begin{equation}
\{x\in\R^{n*} : \ell_j(x)\geq 0 ,\quad j=1,\ldots d\},
\end{equation}
where $\ell_j(x)=\langle u_j,x\rangle-\lambda_j$ for $u_j\in\Z^n$ and
$\lambda_j\in\R$.

Equip $M$ with the Guillemin metric $g$. Then the reduced metric on
$\Delta_0$ is
\begin{equation}
g_{\rm red}= \frac12\sum_{j=1}^d \frac{{d\ell_j}^2}{\ell_j}
\end{equation}
and a dual potential and K{\"a}hler potential are given by
\begin{align}
G&=\frac 12 \sum_{j=1}^d \ell_j(x)\log \ell_j(x)\\
F&=\frac 12 \sum_{j=1}^d \bigl( \lambda_j \log\ell_j(x) +\ell_j(x)\bigr)
\end{align}
\end{thm}
\begin{proof} Using Proposition~\ref{toricquot}, the reduced metric and
dual potential follow immediately from~\eqref{tgred} and~\eqref{tG} by
writing $\mu_j=\ell_j(x)$. It remains to compute the Legendre transform.
The dual coordinates $y$ to the momentum map $x$ are given by
\begin{equation}
y =  \frac{\partial G}{\partial x} = \frac{1}{2}
\sum _{j = 1} ^d u _j \log{\ell _j (x)} + u _j,
\end{equation}
Thus a K{\"a}hler potential is given by
\begin{equation} \begin{split}
- G (x) + \langle x, y \rangle
& = \frac{1}{2} \sum _{j = 1} ^d \bigl(
(\langle u _j, x\rangle  - \ell _j (x)) \log{\ell _j (x)} + \langle
u_j,x\rangle\bigr)\\
& = \frac{1}{2} \sum _{j = 1} ^d \bigl(\ell _j (x)+\lambda_j + \lambda _j 
\log{\ell _j (x)}\bigr).
\end{split} \end{equation}
This differs from $F$ by an additive constant.
\end{proof}

\section*{Appendix: An alternative derivation of the K{\"a}hler potential}

In this Appendix, we show that the expression of the K{\"a}hler
potential $F$ appearing in Theorem \ref{thguillemin} can also be
easily derived from a general expression of the K{\"a}hler potential of
a K{\"a}hler reduction  appearing in \cite{BG}.

We again consider the construction of $M _{\Delta}$ as a K{\"a}hler
reduction of $\mathbb{R} ^{2d} \simeq \mathbb{C} ^d$ with respect to the
standard action of a subgroup $N$ of $T ^d$, cf.~section~\ref{s:toric}.
Notations are the same as in the body of the paper. We denote by $N _{\C}$,
$T ^d _{\C}$, $T ^n _{\C}$ the complexifications of $N$, $T ^d$, $T ^n$; the
complex torus $T ^d _{\C} \simeq \C^* \times \ldots \times \C ^*$ ($d$
times) acts holomorphically on $\R^{2d}\cong \C ^d$ by $(\zeta _1, \cdots,
\zeta _d) \cdot (z _1, \ldots, z _d) = (\zeta _1 z _1, \ldots, \zeta _d z
_d)$, where, for each $j = 1, \ldots, d$, $\zeta _j = t _j e ^{i a _j}$, $t
_j > 0$, is a non-zero complex number; the restriction of this action to the
compact part $T ^d$ coincides with the action already considered. The
restriction of this action to $N _{\C}$ is again a holomorphic extension of
the action of $N$.

Let ${\C} ^d _s$ be the set of points of ${\C} ^d$ whose $N _{\C}$-orbit
cuts $\mu _N ^{-1} (c)$ along a non-empty orbit of $N$ (${\C} ^d _{s}$ is
called the {\it stable} part of ${\C} ^d$ in \cite{VG2}). Then, as a
complex manifold, the K{\"a}hler reduction $M$ coincides with the (ordinary)
quotient of ${\C} ^d _{s}$ by $N _{{\C}}$.

Moreover, for any $z$ in ${\C} ^d _{s}$, there exists a {\it unique}
element $t _z = (t _{z, 1}, \ldots, t _{z, d})$ in ${\rm exp} (i V) \subset
N _{{\C}} \subset T ^d _{{\C}}$ such that $t _z \cdot z$ belongs to $\mu _N
^{-1} (c)$.

We denote by $p$ the natural projection from ${\C} ^d _{s}$ onto $M = {\C}
^d _{s} / N _{{\C}}$; we then have:
\begin{equation}
p (z) = [t _z \cdot z] = [t _{z, 1} z _1, \ldots, t _{z, d} z _d],
\end{equation}
for any $z$ in ${\C} ^d _{s}$, where $[t _z \cdot z]$ denotes the
class mod $N$ of $t _z \cdot z$ in $\mu _N ^{-1} (c)$.

We temporarily assume that the $\lambda _j\in 2\pi \Z$ for all $j$, so that
$\lambda$ can be identified with the character $\chi _{\lambda}$ of $T ^d$
defined by $\chi _{\lambda}\colon [a _1, \ldots, a _d] \to e ^{i
  (\sum _{j = 1} ^d \lambda _j a _j)}$, whereas $- c$ is identified
with the
restriction, $\chi _{- c}$, of $\chi _{\lambda}$ to $N$. We still denote
$\chi _{- c}$ the natural extension of $\chi _{ -c}$ to $N _{{\C}}$, with
values in $\C^*$.

\begin{prop} \label{propA1}
  The pull-back $p ^* \omega$ admits a globally defined $N$ invariant
  potential, $\hat{K}$, on ${\C} ^d _s$, i.e.
\begin{equation} p ^* \omega = d d ^c \hat{K}, \end{equation}
where $\hat{K}$ is given by
\begin{equation} \label{hatK}
\hat{K} (z) = \frac{1}{2} \sum _{j = 1} ^d (\ell _j (x)
+ \lambda _j \log{t _{z, j} ^2}).
\end{equation}
Here, $x = (\nu \circ p) (z)$ is the point of $\Delta$ corresponding to
$p (z)$ and $\ell _j (x) = (x, u _j) - \lambda _j$.
\end{prop}
\begin{proof} According to Theorem 3.1 in \cite{BG}, we have that
\begin{equation}
\hat{K}(z) = K _0 (t _z \cdot z) + \frac{1}{2} \log{|\chi _{-c} (t_z)| ^2},
\end{equation}
where $K _0 = \frac{r ^2}{4}$ is the natural K{\"a}hler potential of the flat
metric of ${\C} ^d _S$. By (\ref{mu}), we then have: $K _0 (t _z \cdot z) =
\frac{1}{4} |t _z \cdot z| ^2 = \frac{1}{2} \sum _{j = 1} ^d \ell _j
(x)$, where $x = (\nu \circ p) (z)$, whereas $|\chi _{-c} (t _z)| ^2 = t
_{z, 1} ^{2 \lambda _1} \ldots t _{z, d} ^{2 \lambda _d}$.
\end{proof}

We define ${\C} ^d _{0}$ as the set of $z$ in ${\C} ^d _s$ such that $z
_j \neq 0$ for all $j$. Its image by $p$ is the open subset, $M _0$, of
$M$ which maps to the interior, $\Delta _0$, of $\Delta$ by $\nu$.

As an easy consequence of Proposition  \ref{propA1}, we get an  alternative
derivation of the Guillemin formula for the K{\"a}hler potential~\cite{VG2}:
\begin{prop} \label{g1} The restriction of $\omega$ to $M _0$ has a
globally defined, $T ^d/N$-invariant potential $K$, i.e.,
\begin{equation} \omega\restr{M _0} = d d ^c K, \end{equation}
with $K = F \circ \nu$, and
\begin{equation} \label{K}
F (x)  = \frac{1}{2} \sum _{j = 1} ^d (\ell _j (x) +
\lambda _j \log{\ell _j (x)}),
\end{equation}
for any $x$ in $\Delta _0$.
\end{prop}
\begin{proof}
  Let $z$ be an element of ${\C} ^d _0$, so that $z _j \neq 0$ for all $j=
  1, \ldots, d$. The element $t _z$, viewed as an element of ${\R} ^+
  \times \ldots \times {\R} ^+$ ($d$ times), is then written as $t _z = (t
  _{z, 1}, \ldots, t _{z, d})$, where, by (\ref{mu}), each positive real
  number $t _{z, j}$ is {\it determined} by: $t _{z, j} ^2 = \frac{2\ell
    _j (x)}{|z _j| ^2}$. We then have:
\begin{equation} \label{t} \begin{split}
\log{\chi _{-c}} (|t _z| ^2)
& = \sum _{j=1} ^d \log{t _{z, j} ^{2 \lambda _j}}
= \sum _{j=1} ^d \lambda _j \log{t _{z, j} ^2} \\
& = \sum _{j=1} ^d \lambda _j (\log{\ell _j} - \log{|z _j| ^2} + \log{2}).
\end{split} \end{equation}
Since the terms in $\log{|z _j| ^2}$ and the constant term $\log{2}$ have
no effect in the computation of $\omega\restr{M _0}$, we directly derive
(\ref{K}) from (\ref{hatK}).
\end{proof}
\begin{rem} {\rm In this approach, we assumed the $\lambda _j$ were
    integral. However, an easy rescaling argument yields the same result
    when the $\lambda _j$ are rational, hence finally, by a continuity
    argument, in the general case where the $\lambda _j$ are real, i.e.,
    for {\it any} Delzant polytope, as we showed in the body of the paper.}
\end{rem}

\end{document}